\documentclass[spanish]{amsart}

\usepackage{amsmath}
\usepackage{amsxtra}
\usepackage{amscd}
\usepackage{amsthm}
\usepackage{amsfonts}
\usepackage{amssymb}
\usepackage{eucal}
\usepackage[all]{xy}
\usepackage{graphicx}
\usepackage{tikz-cd}
\usepackage{mathrsfs}
\usepackage{euler}
\usepackage{hyperref}
\usepackage{color}
\usepackage{longtable}
\usepackage{float}
\usepackage{caption}
\usepackage{enumerate}
\usepackage{faktor}
\usepackage{bbm}
\usepackage{lineno}
\usepackage[english, spanish]{babel}
\usepackage[utf8]{inputenc}

\usepackage[colorinlistoftodos, textsize=tiny]{todonotes}
\def\listtodoname{List of Todos}
\def\listoftodos{\@starttoc{tdo}\listtodoname}

\RequirePackage{color}
\definecolor{myred}{rgb}{0.75,0,0}
\definecolor{mygreen}{rgb}{0,0.5,0}
\definecolor{myblue}{rgb}{0,0,0.65}

\usepackage{tikz}
\usepackage{tikz-cd}
\usetikzlibrary{matrix,arrows,decorations.pathmorphing}

\usepackage{kantlipsum} 

\theoremstyle{plain}
  \newtheorem{theorem}{Teorema}[section]
	
  \newtheorem{proposition}[theorem]{Proposici\'{o}n}
	\newtheorem{prop}[theorem]{Proposici\'{o}n}
	\newtheorem{corollary}[theorem]{Corolario}
	
  \newtheorem{lemma}[theorem]{Lema}

\theoremstyle{definition}
	\newtheorem{definition}[theorem]{Definici\'{o}n}
	
  \newtheorem{example}[theorem]{Ejemplo}

\theoremstyle{remark}
	\newtheorem{remark}[theorem]{Observaci\'{o}n}

\numberwithin{equation}{section}

\DeclareMathAlphabet{\mathcal}{OMS}{cmsy}{m}{n}

\newcommand{\one}{\mathbbm{1}}

\newcommand{\F}{\mathbb{F}}

\renewcommand{\O}{\mathrm{O}}

\newcommand{\Z}{\mathbb{Z}}

\newcommand{\rad}{\mathrm{rad}}
\newcommand{\Char}{\mathrm{car}}
\newcommand{\rank}{\mathrm{rango}}

\newcommand{\Id}{\mathrm{Id}}

\newcommand{\GL}{\mathrm{GL}}

\newcommand{\PO}{\mathrm{PO}}
\newcommand{\PSO}{\mathrm{PSO}}

\newcommand{\SO}{\mathrm{SO}}

\newcommand{\Sp}{\mathrm{Sp}}

\newcommand{\cl}{\mathrm{Cl}}
\newcommand{\arf}{\mathrm{Arf}}
\newcommand{\spin}{\mathrm{sn}}

\newcommand{\paren}[1]{\left(#1\right)}
\newcommand{\set}[1]{\left\{#1\right\}}

\newcommand{\abrac}[1]{\left<#1\right>}
\newcommand{\verts}[1]{\left\lvert#1\right\rvert}

\newcommand\restr[2]{{
  \left.\kern-\nulldelimiterspace 
  #1 
  \right|_{#2} 
  }}

\makeatletter

\DeclareRobustCommand*{\mfaktor}[3][]
{
  { \mathpalette{\mfaktor@impl@}{{#1}{#2}{#3}} }
}
\newcommand*{\mfaktor@impl@}[2]{\mfaktor@impl#1#2}
\newcommand*{\mfaktor@impl}[4]{
  \settoheight{\faktor@zaehlerhoehe}{\ensuremath{#1#2{#3}}}%
  \settoheight{\faktor@nennerhoehe}{\ensuremath{#1#2{#4}}}%
  \raisebox{-0.5\faktor@zaehlerhoehe}{\ensuremath{#1#2{#3}}}%
  \mkern-4mu\diagdown\mkern-5mu%
  \raisebox{0.5\faktor@nennerhoehe}{\ensuremath{#1#2{#4}}}%
}

\newenvironment{abstracts}{%
  \ifx\maketitle\relax
    \ClassWarning{\@classname}{Abstract should precede
      \protect\maketitle\space in AMS document classes; reported}%
  \fi
  \global\setbox\abstractbox=\vtop \bgroup
    \normalfont\Small
    \list{}{\labelwidth\z@
      \leftmargin3pc \rightmargin\leftmargin
      \listparindent\normalparindent \itemindent\z@
      \parsep\z@ \@plus\p@
      
      \itemsep\medskipamount
    }%
}{%
  \endlist\egroup
  \ifx\@setabstract\relax \@setabstracta \fi
}

\newcommand{\abstractin}[1]{%
  \otherlanguage{#1}%
  \item[\hskip\labelsep\scshape\abstractname.]%
}

\makeatother

\title{Grupos ortogonales sobre cuerpos de caracter\'{i}stica positiva}

\author{Robin Zhang}
\address[Robin Zhang]{Department of Mathematics, Columbia University;
Room 509, MC 4406, 2990 Broadway, New York, NY 10027, USA}
\email{rzhang@math.columbia.edu}

\date{16 de mayo de 2022}



\begin{document}

\begin{abstracts}
	\abstractin{spanish}
 	Esta exposici\'{o}n examina la teor\'{i}a de los grupos
	ortogonales y sus subgrupos sobre cuerpos de caracter\'{i}stica
	positiva,
	que recientemente se han utilizado como una herramienta importante en
	el estudio de las formas autom\'{o}rficas y la funcionalidad de
	Langlands.
	Presentamos la clasificaci\'{o}n de grupos ortogonales sobre un
	cuerpo finito $F$
	utilizando la teor\'{i}a de formas bilineales y formas cuadr\'{a}ticas en
	caracter\'{i}stica positiva.
	Usando el determinante y la norma del espinor
	cuando la caracter\'{i}stica de $F$ es impar
	y usando la invariante de Dickson cuando la caracter\'{i}stica de
	$F$ es par, tambi\'{e}n encontramos subgrupos especiales del grupo
	ortogonal.
	
	\abstractin{english}
	This exposition examines the theory of orthogonal groups and
	their subgroups over fields of positive characteristic,
	which has recently been used as an important tool in
	the study of automorphic forms and Langlands functionality.
	We present the classification of orthogonal groups over a finite
	field using the theory of bilinear forms and quadratic forms in
	positive characteristic.
	Using the determinant and spinor norm when the characteristic of
	$F$ is odd and using the Dickson invariant when the characteristic
	of $F$ is even,
	we also look at special subgroups of the orthogonal group.
\end{abstracts}

\selectlanguage{spanish}

\maketitle


\section{Introducci\'{o}n}
La teor\'{i}a de los grupos ortogonales sobre los n\'{u}meros
complejos es muy conocida.
La teor\'{i}a de los grupos ortogonales sobre cuerpos de
caracter\'{i}stica positiva no es tan conocida,
incluso para cuerpos finitos.
Estos grupos siguen siendo importantes en el estudio de la teor\'{i}a
de representaci\'{o}nes, con importantes aplicaciones
recientes en la teor\'{i}a de formas autom\'{o}rficas y la
functorialidad de Langlands
(e.g. \cite{gross-reeder, reeder, lomeli1, lomeli2, castillo}).

Para un espacio vectorial $V$ sobre un cuerpo general $F$,
sea $Q$ una forma cuadr\'{a}tica $Q$ sobre $V$.
En general, el grupo ortogonal de la forma $Q$ es el
grupo de aplicaciones lineales invertibles que conservan $Q$.
Cuando la caracter\'{i}stica de $F$ es impar,
podemos reemplazar la forma cuadr\'{a}tica $Q$ con una forma bilineal
sim\'{e}trica no degenerada $f$ sobre $V$.

Primero, desarrollamos la teor\'{i}a de las formas bilineales para
clasificar los grupos y subgrupos ortogonales sobre un cuerpo
de caracter\'{i}stica impar.
Luego pasamos a la teor\'{i}a de las formas cuadr\'{a}ticas para
clasificar los grupos y subgrupos ortogonales sobre un cuerpo
de caracter\'{i}stica par.
Sobre cuerpos finitos, gran parte del trabajo sobre estos grupos
finitos de tipo Lie fue realizado a principios y mediados del siglo
XX por Dickson, Chevalley, Dieudonn\'{e} y muchos otros.
Dependiendo de si la dimensi\'{o}n de $V$ es par o
impar, hay una o dos clases de equivalencia de tales formas
sobre $V$. Entonces podemos clasificar los tipos de grupos
ortogonales que existen sobre cuerpos de caracter\'{i}stica positiva.
Esta clasificaci\'{o}n est\'{a} contenida esencialmente
en las referencias de Suzuki~\cite[Secci\'{o}n~3.5]{suzuki} y
Wilson~\cite[Cap\'{i}tulo~3]{wilson}.
Uno de los principales objetivos de este art\'{i}culo expositivo es
recopilar estos resultados en el siguiente teorema para beneficio del
lector.

\begin{theorem}[Clasificaci\'{o}n de grupos ortogonales sobre un cuerpo finito]
	\label{thm:classification-orthogonal}
	Sean $F$ un cuerpo finito y $n$ n\'{u}mero entero al menos $2$.
	Los grupos ortogonales sobre $F$ son los siguientes:
	\begin{center}
		\begin{tabular}{ c|c|c| }
			 & $n$ impar & $n$ par \\ 
			\hline
			$\Char(F)$ impar & $\O(n)$ & $\O^\pm(n)$ \\
			\hline
			$\Char(F)$ par & $\Sp(n-1)$ & $\O^\pm(n)$ \\
			\hline
		\end{tabular}
	\end{center}
\end{theorem}

Separamos los casos de Teorema~\ref{thm:classification-orthogonal}
por caracter\'{i}stica del cuerpo $F$.
Repasamos la teor\'{i}a de formas bilineales en
Secci\'{o}n~\ref{subsec:bilinear-forms} y describimos los grupos
ortogonales sobre cuerpos finitos de caracter\'{i}sticas impares en
Secci\'{o}n~\ref{subsec:orthogonal-odd-finite}.
La fila superior del Teorema~\ref{thm:classification-orthogonal} se
prueba en Proposici\'{o}n~\ref{prop:parity-orthogonal-odd}.
De manera similar para cuerpos de caracter\'{i}stica par,
repasamos la teor\'{i}a de formas cuadr\'{a}ticas en
Secci\'{o}n~\ref{subsec:quadratic-forms} y describimos los grupos
ortogonales sobre cuerpos finitos de caracter\'{i}sticas pares en
Secci\'{o}n~\ref{subsec:orthogonal-even-finite}.
La fila inferior del Teorema~\ref{thm:classification-orthogonal} se
prueba en Proposici\'{o}n~\ref{prop:parity-orthogonal-even}.

Tambi\'{e}n encontramos los subgrupos especiales de los grupos
ortogonales.
Cuando $F$ es finito y la caracter\'{i}stica de $F$ es impar,
utilizaremos la norma determinante y la norma espinorial
para definir el subgrupo ortogonal especial y
el n\'{u}cleo espinorial en Secci\'{o}n~\ref{subsec:subgroups-odd}.
Pero cuando la caracter\'{i}stica de $F$ es
par, el determinante es in\'{u}til y no tenemos la norma espinor.
En cambio,
utilizaremos la invariante de Dickson para definir
el subgrupo ortogonal especial y el
n\'{u}cleo espinorial, los cuales coinciden en este caso,
en Secci\'{o}n~\ref{subsec:subgroups-even}.

Algunas de las definiciones y resultados son v\'{a}lidos sobre cuerpos
de caracter\'{i}stica $0$.
Cuando nos enfocamos en cuerpos finitos, muchos de los resultados
tambi\'{e}n se generalizan a cuerpos infinitos.
Sin embargo, los valores num\'{e}ricos y las estructuras de los
grupos ortogonales sobre cuerpos infinitos dependen de las propiedades
del cuerpo, as\'{i} que omitimos esas discusiones.
Para el resto del art\'{i}culo, $F$ es un cuerpo de
caracter\'{i}stica $\Char(F) = p > 0$
y $V$ es un espacio vectorial de
dimensi\'{o}n $n < \infty$ sobre $F$.
Damos pruebas para comodidad del lector,
particularmente cuando usamos ligeras modificaciones
de las referencias citadas.


\section{Clasificac\'{i}on de grupos ortogonales para caracter\'{i}stica impar}

\subsection{Formas bilineales}
\label{subsec:bilinear-forms}
Primero, recordamos algunos hechos sobre las formas bilineales.

\begin{definition}
	Una funci\'{o}n $f: V \times V \rightarrow F$ es
	\emph{una forma bilineal} si y solo si 
	para todos $u, v, w \in V$ y $\lambda \in F$,
	\begin{align*}
		f(\lambda u + v, w) &= \lambda f(u, w) + f(v, w), \\
		f(u, \lambda v + w) &= \lambda f(u, v) + f(u, w).
	\end{align*}
\end{definition}

\begin{definition}
	Sea $f: V \times V \rightarrow F$ una forma bilineal.
	Definimos \emph{el radical} (\emph{izquierdo}) de $f$,
	\[
		\rad(f) := \set{v \in V \mid f(u, v) = 0, \quad \forall u \in V}.
	\]
	
	Adem\'{a}s definimos:
	\begin{itemize}
		\item $f$ es \emph{no singular} cuando $\rad(f) = 0$.
		\item $f$ es \emph{sim\'{e}trica} cuando $f(u, v) = f(v, u)$ para todos $u, v \in V$.
		\item $f$ es \emph{anti-sim\'{e}trica} cuando $f(u, v) = -f(v, u)$ para todos $u, v \in V$.
		\item $f$ es \emph{alternada} cuando $f(v, v) = 0$ para todo $v \in V$.
	\end{itemize}
\end{definition}

\begin{definition}
	Una forma bilineal que es alternada y no singular se llama
	\emph{una forma simpl\'{e}ctica}.
	Un espacio vectorial junto con una forma simpl\'{e}ctica
	se llama \emph{un espacio vectorial simpl\'{e}ctico}.
\end{definition}

\begin{remark}[{cf. \cite[Secci\'{o}n~I.8]{dieudonne},
	\cite[Secci\'{o}n~3.5.6]{suzuki},
	\cite[Secci\'{o}n~3.4.4]{wilson}}]
	\label{rem:alternating-parity}
	Un espacio vectorial
	simpl\'{e}ctico tiene dimensi\'{o}n par y una base
	$\set{e_1, \ldots e_{2m}}$ tal que
	\begin{align*}
			f(e_i, e_{i+m}) &= -f(e_{i+m}, e_i) = 1, \\
			f(e_i, e_j) &= 0 \textrm{ si } i \neq j \pm m.
	\end{align*}
\end{remark}

Los siguientes resultados son conocidos.
\begin{remark}[{\cite[Cap\'{i}tulo~7]{taylor}, \cite[Secci\'{o}n~3.4.1]{wilson}}]
	Sea $f$ una forma bilineal.	
	\begin{enumerate}[(a)]
		\item Si $f$ es alternada, entonces $f$ es anti-sim\'{e}trica.
		\item Si $f$ es anti-sim\'{e}trica y $p \neq 2$, entonces $f$ es alternada.
	\end{enumerate}
\end{remark}

\begin{proof}
	\,
	\begin{enumerate}[(a)]
		\item Tenemos que
			\[
				0 = f(u + v, u + v) = f(u, u) + f(u, v) + f(v, u) + f(v, v) = f(u, v) + f(v, u),
			\]
			entonces $f(u, v) = -f(v, u)$.
		\item Tenemos que $f(v, v) = -f(v, v)$, entonces $2 f(v, v) = 0$. Si $p \neq 2$,
			entonces $f(v, v) = 0$.
	\end{enumerate}
\end{proof}

\subsection{Los grupos ortogonales para \texorpdfstring{$F$}{F} finito}
\label{subsec:orthogonal-odd-finite}

Sea $F$ un cuerpo finito de caracter\'{i}stica impar.
Las dos primeras definiciones tambi\'{e}n son v\'{a}lidas para
otros cuerpos, pero nos centramos en este caso.
Describimos los grupos ortogonales sobre tal cuerpo $F$ y demostramos
la fila superior
de Teorema~\ref{thm:classification-orthogonal}.

\begin{definition}
	Sea $f: V \times V
	\rightarrow F$ una forma bilineal sim\'{e}trica no singular.
	Definimos \emph{el grupo ortogonal},
	\[
		\O(V, f) :=
		\set{g: V \rightarrow V \textrm{ lineal} \mid f\paren{g(u), g(v)} =
		f(u, v), \quad \forall u, v \in V}.
	\]
\end{definition}

\begin{definition}
	Sea $f: V \times V
	\rightarrow F$ una forma simpl\'{e}ctica.
	Definimos \emph{el grupo simpl\'{e}ctico},
	\[
		\Sp(V, f) := \set{g: V \rightarrow V \textrm{ lineal} \mid
			f \paren{g(u), g(v)} = f(u, v) \quad \forall u, v \in V}.
	\]
\end{definition}

\begin{remark}
	Podemos simplemente escribir $\Sp(n) = \Sp(V, f)$ porque existe una
	forma simpl\'{e}ctica can\'{o}nica para $V \cong F^n$.
\end{remark}

Comenzamos con un hecho sobre las formas bilineales simétricas no
singular cuando $F$ es un cuerpo finito $\F_q$ de caracter\'{i}stica
impar.

\begin{lemma}[{cf.~\cite[Secci\'{o}n~3.4.6]{wilson}}]
	\label{lem:equiv-classes-odd}
	Si $F$ es un cuerpo
	finito de caracter\'{i}stica impar,
	entonces hay exactemente dos clases de equivalencia de formas
	bilineales sim\'{e}tricas no singulares
	(bajo la acci\'{o}n) de $\GL_n(F)$).
\end{lemma}

\begin{proof}
	Sea $f$ una forma bilineal sim\'{e}trica sobre $V$.
	Hay vectores $u, v \in V$ tales que $f(u, v) \neq 0$,
	entonces hay un vector $x \in V$ ($x = u$, $x = v$, o $x = u + v$)
	tal que $f(x, x) \neq 0$.
	Escogemos $\alpha \in F$ no cuadrado. Entonces
	\begin{enumerate}
		\item si $f(x, x) = \lambda^2$ es un n\'{u}mero cuadrado,
			$f(x', x') = 1$ para $x' := \lambda^{-1}x$;
		\item si $f(x, x) = \lambda$ no es un n\'{u}mero cuadrado,
			$f(x', x') = \alpha$ para $x' := \sqrt{\frac{\alpha}{\lambda}}x$
			(porque ambos elementos no son residuos cuadr\'{a}ticos).
	\end{enumerate}
	Restringimos $f$ a $(x')^\bot$, continuamos hasta encontrar una base
	perpendicular $\set{e_1, \ldots, e_n}$ tal que
	$f(e_i, e_i) \in \set{0, 1, \alpha}$.
	
	Podemos reducir el n\'{u}mero de $i$ tal que $f(e_i, e_i) = \alpha$
	en dos: si hay $e_i$ y $e_j$ tales que
	$f(e_i, e_i) = f(e_j, e_j) = \alpha$ y $f(e_i, e_j) = 0$,
	podemos elegir $\lambda$ y $\mu$ tales que
	$\lambda^2 + \mu^2 = \frac{1}{\alpha}$ (no cuadrado).
	Entonces $e_i' := \lambda e_i + \mu e_j$ y
	$e_j' := \mu e_i + \lambda e_j$ son una base ortonormal para el
	espacio generado por $e_i$ y $e_j$.
	
	Por lo tanto, siempre podemos encontrar una base $\set{e_1, \ldots e_n}$
	ortogonal tal que
	\begin{enumerate}
		\item $\set{e_1, \ldots, e_n}$ es ortonormal, o
		\item $\set{e_1, \ldots e_{n-1}}$ es ortonormal y
			$f(e_n, e_n) = \alpha$.
	\end{enumerate}
	Entonces hay dos clases de equivalencia.
\end{proof}

Lema~\ref{lem:equiv-classes-odd} implica que hay hasta dos grupos
ortogonales de dimensi\'{o}n fija sobre un cuerpo finito de
caracter\'{i}stica impar.
Cuando hay dos grupos ortogonales distintos, se distinguen por
vectores isotr\'{o}picos.

\begin{definition}
	Un vector $v \in V \backslash \set{0}$ 
	es \emph{un vector isotr\'{o}pico}
	para la forma bilineal $f: V \times V \rightarrow F$
	si $f(v, v) = 0$.
	Un subespacio $W \subset V$ es \emph{totalmente isotr\'{o}pico}
	si $f(u, w) = 0$ para todos $u, w \in W$.
\end{definition}

Por ejemplo (cf. \cite[Secci\'{o}n~3.4.6]{wilson}),
sean $n = 2$ y las dos formas diferentes
$f_1$ y $f_2$ tal que existe una base ortogonal
$\set{x_0, x_1}$ con
\begin{align*}
	f_1(x_0, x_0) &= 1, \\
	f_1(x_1, x_1) &= 1, \\
	f_2(x_0, x_0) &= 1, \\
	f_2(x_1, x_1) &= \alpha \notin (F^\times)^2.
\end{align*}
Si $-1 \in (F^\times)^2$ (i.e.~$\verts{F} \equiv 1 {\pmod 4})$,
entonces
\begin{align*}
	f_1(x_0 + i x_1, x_0 + i x_1) &= 0 \\
	f_2(x_0 + \lambda x_1, x_0 + \lambda x_1) &= 1 + \lambda^2 \alpha \neq 0.
\end{align*}
Si $-1 \notin (F^\times)^2$ (i.e.~$\verts{F} \equiv 3 {\pmod 4})$,
entonces $-\alpha \in (F^\times)^2$ y exista una $\lambda$ tal
que $-\alpha = \lambda^{-2}$, por lo que entonces
\begin{align*}
	f_1(x_0 + i x_1, x_0 + i x_1) &\neq 0 \\
	f_2(x_0 + \lambda x_1, x_0 + \lambda x_1) &= 1 + \lambda^2 \alpha = 0.
\end{align*}
En ambos casos, exactamente una de las dos $f_i$ tiene un vector
isotr\'{o}pico.
Esta es una distinci\'{o}n geom\'{e}trica que diferencia los
dos grupos ortogonales.

Para $n$ par en general, 
exactamente una de las dos $f_i$ tiene un subespacio
totalmente isotr\'{o}pico de dimension $\frac{n}{2}$,
mientras que el otra tiene un subespacio totalmente
isotr\'{o}pico de dimension $\frac{n}{2} - 1$.
Escribimos $\O^+(n)$ para $\O(V, f)$ si
existe un subespacio totalmente isotr\'{o}pico de dimension
$\frac{n}{2}$, y $\O^-(n)$ de lo contrario.
			
Para completar la prueba de la fila superior
de Teorema~\ref{thm:classification-orthogonal},
demostramos que tenemos exactamente uno o exactamente dos grupos
ortogonales diferentes
dependiendo de la paridad de la dimensi\'{o}n $n$.
\begin{prop}[{cf.~\cite[Secci\'{o}n~3.5.8]{suzuki}}]
	\label{prop:parity-orthogonal-odd}
	\,
	\begin{enumerate}
		\item Si $n$ es impar, entonces hay solo un grupo ortogonal
			(excepto por isomorfismo) $\O(n)$
		\item Si $n$ es par, entonces hay dos grupos 
			(excepto por isomorfismo) $\O^\pm(n)$
	\end{enumerate}
\end{prop}

\begin{proof}
	\,
	\begin{enumerate}
		\item Sean $f$ y $f'$ dos formas bilineales sim\'{e}tricas no
			singulares. Escribimos $n = 2m + 1$.
			Usando un proceso similar al de la prueba de
			Lema~\ref{lem:equiv-classes-odd}, hay
			$\alpha, \beta \in F^\times$ y bases
			$\set{e_1, \ldots, e_{2m+1}},
			\set{e_1', \ldots, e_{2m+1}'} \subset V$
			tales que
			\begin{align*}
				f(e_{2i-1}, e_{2i}) &= 1 & \textrm{ si } 1 \leq i \leq m \\
				f(e_{2m+1}, e_{2m+1}) &= \alpha \\
				f(e_i, e_j) &= 0 & \textrm{ para otros } (i, j), \\
				f'(e_{2i-1}', e_{2i}') &= 1 & \textrm{ si } 1 \leq i \leq m \\
				f'(e_{2m+1}', e_{2m+1}') &= \beta \\
				f'(e_i', e_j') &= 0 & \textrm{ otros } (i, j).
			\end{align*}
			Definimos $\sigma \in \GL(V)$ por
			\begin{align*}
				e_{2i-1} &\mapsto e_{2i-1}', &\textrm{ si } 1 \leq i \leq m), \\
				e_{2i} &\mapsto \alpha^{-1} \beta e_{2i}' &\textrm{ si } 1 \leq i \leq m), \\
				e_{2m+1} &\mapsto e_{2m+1}'.
			\end{align*}
			Entonces $f' = \alpha^{-1} \beta f^\sigma$, con
			$f^\sigma(u, v) := f\paren{\sigma^{-1}(u), \sigma^{-1}(v)}$.
			La multiplicaci\'{o}n por una constante y la acci\'{o}n de $\GL(V)$ no
			cambia el grupo ortogonal
			(cf. \cite[Secci\'{o}n~3.5.3]{suzuki}), i.e.
			$\O(V, f') = \O(V, \alpha^{-1} \beta f^\sigma) = \O(V, f)$.
		\item Si $n = 2m$,
			hay una base
			$\set{e_1, \ldots, e_{2m}} \subset V$
			tal que
			\begin{align*}
				f(e_{2i-1}, e_{2i}) &= 1 & \textrm{ si } 1 \leq i \leq m-1 \\
				f(e_i, e_j) &= 1 & \textrm{ para otros } (i, j).
			\end{align*}
			Hay dos casos, que corresponden a a la existencia de un
			subespacio totalmente isotr\'{o}pico de dimension $m$:
			\begin{enumerate}
				\item $f(e_{2m-1}, e_{2m-1}) = f(e_{2m}, e_{2m}) = 0$ y
					$f(e_{2m-1}, e_{2m})=1$, o
				\item $f(e_{2m-1}, e_{2m-1}) = -\alpha f(e_{2m}, e_{2m})
					\neq 0$ y $f(e_{2m-1}, e_{2m}) = 0$.
			\end{enumerate}
			
			En el caso (a), hay un
			subespacio totalmente isotr\'{o}pico de dimension $m$
			y el orden del grupo ortogonal $\O(V, f) = \O^+(n)$ es
			\[
				\verts{\O^+(n)} = 2 (q^m - 1)q^{m(m-1)}\prod_{i=1}^{m-1}(q^{2m-2i} - 1).
			\]
			En el caso (b), no hay un
			subespacio totalmente isotr\'{o}pico de dimension $m$
			y el orden del grupo ortogonal $\O(V, f) = \O^-(n)$ es
			\[
				\verts{\O^-(n)} = 2 (q^m + 1)q^{m(m-1)}\prod_{i=1}^{m-1}(q^{2m-2i} - 1).
			\]
			Los dos grupos claramente no son isomorfos porque 
			sus \'{o}rdenes son diferentes.
	\end{enumerate}
\end{proof}



\section{Clasificac\'{i}on de grupos ortogonales para caracter\'{i}stica \texorpdfstring{$2$}{2}}

\subsection{Formas cuadr\'{a}ticas}
\label{subsec:quadratic-forms}
De manera similar con las formas bilineales,
podemos definir grupos ortogonales en t\'{e}rminos de
formas cuadr\'{a}ticas. Para cuerpos de caracter\'{i}stica impar,
esto es equivalente a la formulaci\'{o}n que utiliza formas bilineales
sim\'{e}tricas.
Cuando trabajemos con cuerpos de caracter\'{i}stica $2$, solo usaremos
formas cuadr\'{a}ticas.

Sean $F$ un cuerpo de caracter\'{i}stica $p$
(no necesariamente par) y
$V$ un espacio vectorial sobre $F$ de dimensi\'{o}n $n$.

\begin{definition}
	Una funci\'{o}n $Q: V \rightarrow F$ es
	\emph{una forma cuadr\'{a}tica} cuando 
	para todos $u, v \in V$ y $\lambda \in F$,
	\begin{align*}
		Q(\lambda u + v) &= \lambda^2 Q(u) + \lambda f_Q (u, v) + Q(v),
	\end{align*}
	
	\noindent
	donde $f_Q: V \times V \rightarrow F$ es una forma bilineal
	sim\'{e}trica.
\end{definition}

Una forma cuadr\'{a}tica $Q: V \rightarrow F$
determina una forma bilineal sim\'{e}trica $f_Q$.
De hecho, $f_Q$ es alternada cuando $p = 2$
porque
\[
	0 = Q(2v) = 2Q(v) + f_Q(v, v) = f_Q(v, v).
\]
Dada una forma bilineal sim\'{e}trica $f: V \times V \rightarrow F$,
también podemos obtener una forma cuadr\'{a}tica
$Q_f(v) := \frac{f(v, v)}{2}$, pero solo cuando $p \neq 2$.
Por lo tanto, solo existe una equivalencia entre la teor\'{i}a de formas
bilineales sim\'{e}tricas y la teor\'{i}a de formas cuadr\'{a}ticas cuando
la caracterist\'{i}ca de $F$ no es $2$.

\begin{definition}
	Sea $Q: V \rightarrow F$ una forma cuadr\'{a}tica.
	Definimos \emph{el radical} de $Q$,
	\[
		\rad(Q) := \set{v \in \rad(f_Q) \mid Q(v) = 0},
	\]
	y \emph{el \'{a}lgebra de Clifford},
	\[
		\cl(V, Q) := \bigoplus_{n \geq 0} V^{\otimes n} / \abrac{v \otimes v - Q(v) \otimes \one}.
	\]
	
	Adem\'{a}s definimos:
	\begin{itemize}
		\item $Q$ es \emph{no singular} cuando $\rad(Q) = 0$.
		\item $Q$ es \emph{no degenerada} cuando $f_Q$ es no singular.
	\end{itemize}
\end{definition}

Introducimos la invariante
de Arf para describir las formas cuadr\'{a}ticas sobre cuerpos
de caracter\'{i}stica par.
Mencionamos la invariante de Arf porque las formas cuadr\'{a}ticas no
singulares sobre un espacio vectorial sobre $F$ est\'{a}n
determinadas por la invariante de Arf y el \'{a}lgebra de Clifford.
Para el resto de la secci\'{o}n, $F$ es un cuerpo de
caracter\'{i}stica par

Sea $Q$ una forma cuadr\'{a}tica no degenerada (entonces $n$ es par).
Definimos la relaci\'{o}n de equivalencia $\sim$:
$Q \sim Q'$ cuando existe $C \in \GL_n(F)$ tal que
$Q(v) = Q'(Cv)$.
Por un teorema de Arf \cite[Satz~2]{arf}
(cf. \cite[Secci\'{o}n~4]{lorenz-roquette}),
existen formas cuadr\'{a}ticas binarias $Q_i$ tales que
\[
	Q \sim \bigoplus_i^r Q_i.
\]

\begin{definition}[{cf. \cite[Secci\'{o}n~7]{lorenz-roquette}}]
	Sea $U$ el subgrupo aditivo de $F$ definido por
	\[
		U := \set{u^2 + u \mid u \in F}.
	\]
	Definimos \emph{la invariante de Arf} 
	por
	\[
		\arf(ax^2 + xy + by^2) := ab + U \in F/U.
	\]
	Para una forma cuadr\'{a}tica no degenerada $Q$,
	definimos \emph{la invariante de Arf}
	(o \emph{el pseudodeterminante} \cite{dieudonne-1955}) de $Q$ por
	\[
		\arf(Q) := \sum_{i=1}^r \arf(Q_i).
	\]
\end{definition}

\begin{remark}
	Para cualquier forma cuadr\'{a}tica binaria $Q$, existen
	$a$ y $b$ en $F$ tales que $Q(x, y) \sim ax^2 + xy + by^2$.
	La invariante de Arf est\'{a} bien definido porque
	si $ax^2 + xy + by^2 \sim a'x^2 + xy + b'y^2$, existe
	$u \in F$ tal que $ab - a'b' = u^2 + u$.
\end{remark}

\begin{theorem}[{\cite[Teoremas~11-12]{arf}, \cite{baeza}}]
	Si $F$ es perfecto, una forma cuadr\'{a}tica no singular
	$Q$ est\'{a} determinada \'{u}nicamente por el par
	$(\arf(Q), n)$.
	
	Si $F$ no es perfecto y $[F:F^2] \leq 2$, una forma cuadr\'{a}tica
	no singular	$Q$ est\'{a} determinada \'{u}nicamente por el triplete
	$(\arf(Q), n, \cl(V, Q))$ donde $\cl(V, Q)$ es el \'{a}lgebra
	de Clifford.
\end{theorem}

\begin{remark}
	Para una excelente revisi\'{o}n de la historia de los teoremas de Arf y
	la correcci\'{o}n de sus errores,
	consulte \cite[Secci\'{o}n~9]{lorenz-roquette}.
\end{remark}

\begin{example}
	Si $k$ es un cuerpo perfecto, entonces $F = k(t)$ y $F = k((t))$ son
	cuerpos tales que $F$ no es perfecto y $[F:F^2] \leq 2$.
\end{example}

\subsection{Los grupos ortogonales para \texorpdfstring{$F$}{F} finito}
\label{subsec:orthogonal-even-finite}

Sea $F$ un cuerpo finito de caracter\'{i}stica $2$.
La primera definicion tambi\'{e}n es v\'{a}lida para
otros cuerpos, pero nos centramos en este caso.
Describimos los grupos ortogonales sobre tal cuerpo $F$ y demostramos
la fila inferior
de Teorema~\ref{thm:classification-orthogonal}.

\begin{definition}
	Sea $Q: V \rightarrow F$ una forma cuadr\'{a}tica no singular.
	Definimos \emph{el grupo ortogonal} de $V$ sobre $Q$ como
	\[
		\O(V, Q) :=
		\set{g: V \rightarrow V \textrm{ lineal} \mid Q\paren{g(v)} =
		Q(v), \quad \forall v \in V}.
	\]
\end{definition}

Los siguientes resultados son conocidos
(cf. \cite[Secci\'{o}n~3.4.7 y Secci\'{o}n~3.8]{wilson}).
\begin{lemma}
	\label{lem:multiplicity-one}
	Si $Q: V \rightarrow F$ es una forma cuadr\'{a}tica degenerada y no
	singular,
	entonces $\dim \rad(f) = 1$.
\end{lemma}

\begin{proof}
	La forma bilineal sim\'{e}trica $f_Q$ es alternada porque
	$\Char(F) = 2$. Si $u, v \in \rad(f)$,
	\[
		Q(\lambda u + v) = \lambda^2 Q(v) + Q(w),
	\]
	entonces $\restr{Q}{\rad(f)}$ es una aplicaci\'{o}n semilineal
	a $F$ y $\rad(Q)$ es un subespacio de codimensi\'{o}n
	$0$ o $1$ de $\rad(f)$. Pero $\rad(f) \neq 0$ y $\rad(Q) = 0$
	porque $Q$ es degenerada y no singular, entonces
	$\dim \rad(f) = 1$.
\end{proof}

\begin{lemma}
	\label{lem:lift}
	Si $Q: V \rightarrow F$ es una forma cuadr\'{a}tica degenerada
	no singular, entonces
	$f_Q$ induce una forma simpl\'{e}ctica sobre $V/\rad(f_Q)$
	tal que
	$\O(V, Q) \cong \O(V /\rad(f_Q), f_Q')$.
\end{lemma}

\begin{proof}
	Para el complemento $W \cong V/\rad(f_Q)$ de $V/\rad(f_Q)$ en $V$,
	$V = W \oplus \rad(f_Q)$. Para $g \in \O(V, Q)$ y
	$w, w' \in W$,
	\begin{align*}
		Q(w) + f_Q(w, w') + Q(w') &= Q(w + w') \\
			&= Q(\restr{g}{W}(w + w')) \\
			&= Q(\restr{g}{W}(w)) + f_Q(\restr{g}{W}(w), \restr{g}{W}(w'))
				+ Q(\restr{g}{W}(w')) \\
			&= Q(w) + f_Q(\restr{g}{W}(w), \restr{g}{W}(w')) + Q(w').
	\end{align*}
	Entonces, $f_Q(w, w') = f_Q(\restr{g}{W}(w), \restr{g}{W}(w'))$
	para todos $w, w' \in W$. Usando $\restr{g}{W}$ y el isomorfismo
	$W \cong V/\rad(f_Q)$,	
	hay $g' \in \O(V /\rad(f_Q), f_Q')$.
	
	El elemento de $\O(V /\rad(f_Q), f_Q')$ obtenido de cada
	$g \in \O(V, Q)$ es \'{u}nico. Sean $g_1, g_2 \in
	\O(V, Q)$ tales que $g_1' = g_2' \in \O(V /\rad(f_Q), f_Q')$.
	La dimensi\'{o}n de
	$\rad(f_Q)$ es $1$, entonces podemos escribir
	$\rad(f_Q) = \abrac{u}$.
	Sean $v \in V$ y $w \in W$.
	\begin{align*}
		f_Q(v, g_1(w) + g_2(w)) &= f_Q(v, g_1'(w) + g_2'(w)) \\
			&= f_Q(v, 0) = 0. \\
		Q(g_1(w) + g_2(w)) &= Q(g_1(w)) + f_Q(g_1(w), g_2(w)) + Q(g_2(w)) \\
			&= 2Q(w) + f_Q(g_1'(w), g_2'(w)) \\
			&= f_Q(g_1'(w), g_1'(w)) \\
			&= 0.
	\end{align*}
	Entonces, $g_1(w) + g_2(w) \in \rad(Q) = \set{0}$ y
	$g_1(w) = g_2(w)$ para todos $w \in W$.
	Solo necesitos verificar que
	$\restr{g_1}{\abrac{u}} = \restr{g_2}{\abrac{u}}$
	porque $\restr{g_1}{W} = \restr{g_2}{W}$. Pero,
	$g_i(u) \in \abrac{u}$ porque $g_i$ es lineal.
	Adem\'{a}s,
	\begin{align*}
		Q(g_i(u) + u)) &= Q(g_i(u)) + f_Q(g_i(u), u) + Q(u) \\
			&= 2Q(u) + f_Q(g_i(u), u) \\
			&= 0.
	\end{align*}
	Entonces $g_i(u) + u \in \rad(Q) = \set{0}$ y
	$g_i(u) = u$ para cada $i \in \set{1, 2}$.
\end{proof}

\begin{proposition}
\label{fact:char2-symplectic}
	Si $Q: V \rightarrow F$ es una forma cuadr\'{a}tica
	no singular, entonces
	\begin{align*}
		\O(V, Q) &\cong \Sp(n-1) &\textrm{si } n \textrm{ es impar}, \\
		\O(V, Q) &\leq \Sp(n) &\textrm{si } n \textrm{ es par}.
	\end{align*}
\end{proposition}

\begin{proof}
	Sea $g \in \O(V, Q)$ y $u, v \in V$.
	Tenemos las igualdades,
	\begin{align*}
		Q(u) + f_Q(u, v) + Q(v) &= Q(u + v) \\
			&= Q(g(u + v)) \\
			&= Q(g(u) + g(v)) \\
			&= Q(g(u)) + f_Q(g(u), g(v)) + Q(g(v)) \\
			&= Q(u) + f_Q(g(u), g(v)) + Q(v).
	\end{align*}
	Entonces, $f_Q(u, v) = f_Q(g(u), g(v))$ para todos $u, v \in V$
	y $g \in \O(V, f_Q)$. Por lo tanto, $\O(V, Q) \leq \Sp(n)$. 
	
	Si $n$ es impar, $\rad(f_Q) \neq 0$
	(de lo contrario $f_Q$ es simpl\'{e}ctica).
	Por Lema~\ref{lem:multiplicity-one}, $\dim \rad(f_Q) = 1$.
	Por Lema~\ref{lem:lift}, $\O(V, Q) \cong \O(V/\rad(f_Q), f_Q)$.
	La forma bilineal alternada $f_Q$ induce una forma simpl\'{e}ctica
	sobre $V/\rad(f_Q)$, entonces
	\[
		\O(V/rad(f_Q), f_Q') \cong \Sp(n-1).
	\]
\end{proof}

Proposici\'{o}n~\ref{fact:char2-symplectic} da la esquina inferior
izquierda de Teorema~\ref{thm:classification-orthogonal}.
Para completar la demostraci\'{o}n de la fila inferior
de Teorema~\ref{thm:classification-orthogonal},
solo necesitamos contar el n\'{u}mero de grupos ortogonales de
dimensi\'{o}n par.
Como en Secci\'{o}n~\ref{subsec:orthogonal-odd-finite},
usamos vectores isotr\'{o}picos.

\begin{definition}
	Un vector $v \in V \backslash \set{0}$ 
	es \emph{un vector isotr\'{o}pico}
	para la forma cuadr\'{a}tica $Q: V \rightarrow F$
	si $Q(v) = 0$.
\end{definition}

\begin{lemma}[{cf.~\cite[Secci\'{o}n~3.5.9]{suzuki}, \cite[Secci\'{o}n~3.4.7]{wilson}}]
	\label{lem:equiv-classes-even}
	Hay exactemente dos clases de equivalencia de formas
	cuadr\'{a}ticas no singulares.
\end{lemma}

\begin{proof}
	Construimos una base de la misma manera que 
	la base simpl\'{e}ctica en
	Observaci\'{o}n~\ref{rem:alternating-parity}, pero
	con $Q(e_i) = Q(f_i) = 0$ cuando sea posible.
	Si $n > 2$ y $Q(e_i) \neq 0$, podemos reemplazar $e_i$
	par
	\[
		e_i' = \paren{\frac{Q(f_i)}{Q(e_i)}}^\frac{\verts{F}}{2} e_i + f_i,
	\]
	si $Q(e_i) \neq 0$ (en cuyo caso $Q(e_i') = 0$). Si $Q(e_i) = 0$,
	podemos reemplazar $f_i$ par $f_i' = f_i + Q(f_i)e_i$ para
	garantizar que $Q(f_i') = 2Q(f_i) = 0$.
	
	Pero para los dos primeros elementos de la base, podemos
	elegir $e_1$ y $f_1$ tal que $Q(e_1) = f(e_1, f_1) = 1$,
	$f(e_1, f_1 + \lambda e_1) = 1$ y
	$Q(f_1 + \lambda e_1) = Q(f_1) + \lambda^2 + \lambda$.
	Podemos reemplazar $f_1$ par $f_1' = f_1 + \lambda e_1$.
	Hay exactamente dos clases de equivalencia de formas
	cuadr\'{a}ticas no singulares, dependiendo de si
	$Q$ tiene un vector isotr\'{o}pico (i.e.
	si $Q(f_1') = \lambda^2 + \lambda + Q(f_1) = 0$ tiene una
	soluci\'{o}n en $F$).	
\end{proof}

Como en Proposici\'{o}n~\ref{prop:parity-orthogonal-odd},
las dos clases de equivalencia de formas cuadr\'{a}ticas no singulares
en Lema~\ref{lem:equiv-classes-even} producen dos clases de
isomorfismo de grupos ortogonales con
diferentes ordenes.
\begin{proposition}[{cf.~\cite[Secci\'{o}n~3.5.10]{suzuki}}]
	\label{prop:parity-orthogonal-even}
	\,
	\begin{enumerate}
		\item Si $n$ es impar, entonces hay solo un grupo ortogonal
			(excepto por isomorfismo) $\O(n) \cong \Sp(n-1)$
		\item Si $n$ es par, entonces hay dos grupos 
			(excepto por isomorfismo) $\O^\pm(n)$
	\end{enumerate}	
\end{proposition}

\begin{proof}
	\,
	\begin{enumerate}
		\item Por Proposici\'{o}n~\ref{fact:char2-symplectic}.
		\item 
		Las dos clases de equivalencia de formas cuadr\'{a}ticas no singulares
		en Lema~\ref{lem:equiv-classes-even} se distinguen por diferentes
		n\'{u}meros de vectores isotr\'{o}picos.
		Como en Proposici\'{o}n~\ref{prop:parity-orthogonal-odd},
		el orden del grupo ortogonal $\O(V, Q) = \O^+(n)$ en un caso es
		\[
			\verts{\O^+(n)} = 2 (q^m - 1)q^{m(m-1)}\prod_{i=1}^{m-1}(q^{2m-2i} - 1).
		\]
		En el otro caso, el orden del grupo ortogonal $\O(V, Q) = \O^-(n)$
		es
		\[
			\verts{\O^-(n)} = 2 (q^m + 1)q^{m(m-1)}\prod_{i=1}^{m-1}(q^{2m-2i} - 1).
		\]
	\end{enumerate}
\end{proof}


\section{Subgrupos de grupos ortogonales}

\subsection{Los subgrupos para caracter\'{i}stica impar}
\label{subsec:subgroups-odd}

Sea $F$ un cuerpo finito de caracter\'{i}stica impar.
Todos los elementos del grupo ortogonal tienen determinante $\pm 1$.
Por lo tanto tenemos \emph{el grupo ortogonal especial}
$\SO^{(\pm)}(n)$, el subgrupo de $\O^{(\pm)}(n)$ de
\'{i}ndice $2$ definido como el n\'{u}cleo del determinante.
Al igual que con el grupo ortogonal en
Teorema~\ref{thm:classification-orthogonal}, hay una clase de
isomorfismo de grupos ortogonales especiales $\SO(n)$ si $n$ es impar
y hay dos clases de isomorfismo de grupos ortogonales especiales
$\SO^{\pm}$ si $n$ es par.

Tambi\'{e}n $\pm \Id \in \O^{(\pm)}(n)$
(con $-\Id \in \SO^{(\pm)}(n)$ si y solo si $n$ es par),
entonces tenemos los grupos de cocientes $\PO^{(\pm)}(n)
:= \O^{(\pm)}(n) / \set{\pm \Id}$ y
$\PSO^{(\pm)}(n) := \SO^{(\pm)}(n) / \set{\pm \Id}$
llamados respectivamente
\emph{el grupo ortogonal proyectivo} y
\emph{el grupo ortogonal especial proyectivo}.

\begin{definition}
	\label{def:spinor}
	Para una forma bilineal sim\'{e}trica $f: V \times V \rightarrow F$ y
	un vector $v \in V$ tal que $f(v, v) \neq 0$,
	definimos la reflexi\'{o}n
	\[
		r_v: x \mapsto x - 2\frac{f(x, v)}{f(v, v)} v.
	\]
	Un elemento $g$ del $\SO^{(\pm)}(n)$ tiene \emph{la norma espinor} $1$
	cuando $g$ es un producto de un n\'{u}mero par
	de reflexiones de vectores de norma $1$
	y un n\'{u}mero par de reflexiones de vectores de norma
	no cuadrada.
	Un elemento $g$ del $\SO^{(\pm)}(n)$ tiene la norma espinor $-1$
	de lo contrario.
\end{definition}

Usando el n\'{u}cleo de la norma espinor
$\spin: \SO^{(\pm)}(n) \rightarrow \set{\pm 1}$,
definimos los subgrupos de \'{i}ndice $2$
\begin{align*}
	\Omega^{(\pm)}(n) :&=
		\ker(\spin) \leq \SO^{(\pm)}(n), \\
	\mathrm{P}\Omega^{(\pm)}(n) :&=
		\ker\paren{\restr{\spin}{\PSO^{(\pm)}}} \leq \PSO^{(\pm)}(n).
\end{align*}

\subsection{Los subgrupos para caracter\'{i}stica \texorpdfstring{$2$}{2}}
\label{subsec:subgroups-even}

En esta secci\'{o}n, $F$ es cuerpo de caracter\'{i}stica $2$.
Podemos considerar que $F$ es finito, pero las siguientes definiciones
y resultados tambi\'{e}n son v\'{a}lidos para cuerpos infinitos.

Para caracter\'{i}stica $2$, el determinante es in\'{u}til y
no tenemos la norma espinor (cf. Definici\'{o}n~\ref{def:spinor}).
Sea $Q: V \rightarrow F$ una forma cuadr\'{a}tica no singular.
Definimos la transvecci\'{o}n ortogonal
(\guillemotleft la reflexi\'{o}n \guillemotright)
\[
	t_v: x \mapsto x + \frac{f_Q(x, v)}{f_Q(v, v)} v
\]
para $v \in V$ tal que $f_Q(v, v) \neq 0$.
El transvecci\'{o}n $t_v$ es lineal y conserva $Q$ porque
\[
	Q\paren{x + f(x, v)v} = Q(x) + f_Q(x, v)^2 + f_Q(x, v)^2 Q(v) = Q(x).
\]

\begin{theorem}[Cartan--Dieudonn\'{e}--Kneser]
	El grupo ortogonal $\O(V, Q)$ es generado por $\set{t_v}_{v \in V}$,
	excepto cuando $F = \F_2$ y $n = 4$.
\end{theorem}

\begin{remark}
	Wilson~\cite[Cap\'{i}tulo~3]{wilson} no menciona la \'{u}nica
	excepci\'{o}n dada por Cartan--Dieudonn\'{e}--Kneser (cf.
	\cite[Teorema~I.5.1]{chevalley}, \cite[pp.20--22]{dieudonne},
	\cite{kneser2}, \cite[Teorema~6.3.4]{knus}, y
	\cite[Teorema~11.42]{taylor}).
\end{remark}

\begin{definition}
	Definimos \emph{la invariante de Dickson} para
	\[
		g = \prod_{i=1}^r t_{v_i},
	\]
	como
	\[
		D(g) :=
		\begin{cases}
			1 & r \textrm{ es impar,} \\
			0 & r \textrm{ es par.}
		\end{cases}
	\]
	Equivalente, $D(g) = \rank\paren{\Id - g} (\textrm{m\'{o}dulo } 2)$.
\end{definition}

La invariante de Dickson tambi\'{e}n se le conoce como
\emph{el cuasideterminante}	o \emph{el pseudodeterminante}
porque son esencialmente el mismo que el determinante
en caracter\'{i}stica impar.
\begin{theorem}[{\cite[Teorema~2]{dye}, cf. \cite[Teorema~11.43]{taylor}}]
	\label{thm:dye}
	La invariante de Dickson es un homomorfismo
	\[
		D: \O(V, Q) \rightarrow \Z/2\Z.
	\]

	Adem\'{a}s para caracter\'{i}stica diferente a $2$,
	\[
		\det(g) = (-1)^{D(g)}.
	\]
	Entonces en este caso,
	el determinante es equivalente a la invariante de Dickson.
\end{theorem}

Podemos usar la invariante de Dickson como determinante para definir
subgrupos.

\begin{definition}
	Definimos \emph{el grupo ortogonal especial} como el n\'{u}cleo de $D$,
	\[
		\SO(V, Q) := \set{g \in \O(V, Q) \mid D(g) = 0}.
	\]
	Este grupo también se conoce como \emph{el	n\'{u}cleo espinorial}
	$\Omega(V, Q)$.
\end{definition}

Una consecuencia de la segunda parte del Teorema~\ref{thm:dye},
que en realidad funciona para un cuerpo general $F$,
explica el nombre \guillemotleft grupo ortogonal especial \guillemotright.
\begin{corollary}
	Si $F$ es un cuerpo de caracter\'{i}stica $0$ o impar,
	el n\'{u}cleo de $D$ es isomorfo al subgrupo de elementos
	de determinantes $1$ en $\O(V, Q)$ y tambi\'{e}n al subgrupo
	de elementos de norma espinorial $1$ en $\O(V, Q)$.
\end{corollary}

\begin{theorem}[\cite{kneser, pollak1962, connors, pollak1970}]
	\label{thm:commutator}
	El subgrupo commutador
	\[
		\O(V, Q)^{(1)} := [\O(V, Q), \O(V, Q)]
	\]
	es isomorfo al $\Omega(V, Q)$.
\end{theorem}

\begin{remark}
	Si $F$ es un cuerpo global o local, la afirmaci\'{o}n del
	Teorema~\ref{thm:commutator} sigue siendo cierto
	con una sola excepci\'{o}n: si $F$ es un cuerpo local
	de caracter\'{i}stica impar, $n = 4$, y
	$V$ es anisotr\'{o}pico. En este caso,
	Pollak~\cite{pollak1962} demostr\'{o} que
	$[\Omega(V, Q):\O(V, Q)^{(1)}] = 2$.
\end{remark}


\section{Agradecimientos}
Estas notas se basan en una charla impartida en el taller
de Casa Matem\'{a}tica Oaxaca--Banff International Research Station
(CMO--BIRS),
\guillemotleft Teor\'{i}a de n\'{u}meros en Am\'{e}rica \guillemotright,
en agosto de 2019.
Agradecemos su apoyo, as\'{i} como las discusiones con Lea Beneish,
Michael Harris, Luis Lomel\'{i} y Alberto M\'{i}nguez.
Tambi\'{e}n agradecemos a los \'{a}rbitros an\'{o}nimos por sus
\'{u}tiles comentarios.

El autor fue apoyado por la
National Science Foundation Graduate Research
Fellowship Program, Grant No. DGE-1644869. Todas las opiniones,
hallazgos y conclusiones o recomendaciones expresadas en este material
pertenecen al autor y no reflejan necesariamente los puntos de vista
de la National Science Foundation.


\nocite{*}
\bibliography{bibliography}{}
\bibliographystyle{amsalpha}

\end{document}